\documentstyle[amssymb,amsmath,eucal]{article}

\begin{document}

\def\R {{\Bbb R }}
\def\C {{\Bbb C }}

\def\fg{{\frak g}}
\def\fp{{\frak p}}
\def\fn{{\frak n}}
\def\fr{{\frak r}}
\def\fq{{\frak q}}
\def\fh{{\frak h}}
\def\fm{{\frak m}}
\def\Ad{{\rm Ad}}

\begin{center}
{\large\bf
A construction of
finite-dimensional faithful representation
of Lie algebra}

\smallskip

\sc Yurii A.Neretin

\end{center}

The Ado theorem is a fundamental fact,
which has a reputation
 to be a 'strange theorem'.
 We give its natural proof.

\smallskip

{\bf 1. Construction of faithful representation.}
Consider a finite-dimensional
 Lie algebra $\fg$. Assume that
 $\fg$ is a semidirect product $\fp\ltimes\fn$
 of a subalgebra $\fp$ and a nilpotent
 ideal $\fn$.
Assume that the adjoint action of $\fp$ on $\fn$
is faithful, i.e., for any $z\in\fp$,
there exists $x\in\fn$ such that
$[z,x]\ne 0$.

Consider the minimal $k$
such that all the commutators
$$[\dots[[x_1,x_2],x_3],\dots,x_k], \qquad
x_j\in\fn$$
are 0.

Denote by $U(\fn)$ the enveloping
algebra of $\fn$.
The algebra $\fn$ acts on $U(\fn)$
by the left multiplications.
The algebra $\fp$ acts on $U(\fn)$
by the derivations
$$
d_z x_1x_2 x_3\dots x_l=
[z,x_1]x_2x_3\dots x_l+x_1[z,x_2]x_3\dots x_l+\dots,
\qquad
\text{where $z\in \fp$}
.$$
This defines the action of the semidirect product
$\fp\ltimes\fn=\fg$ on $U(\fn)$.

 Denote by $I$
the subspace in $U(\fn)$
spanned by all the products
$x_1 x_2\dots x_N$, where $N>k+2$.
Obviously,

 1. $I$ is the two-side ideal in $U(\fn)$.

 2.
Consider the linear span $\cal A$
 of all the elements having the form
 $1$, $x$, $x_1x_2\in U(\fn)$.
 Obviously, $I\cap {\cal A}=0$.

 3. $I$ is invariant with respect to the
 derivations $d_z$.

Obviously, the
module $U(\fn)/I$
is a finite-dimensional faithful module
over $\fg$.

\smallskip

{\bf 2. The Ado theorem.}

\smallskip

{\sc Lemma 1.} {\it
Any finite-dimensional
 Lie algebra $\fq$
admits an embedding to an algebra $\fg$
such that

 a) $\fg$
is a semidirect product
of    a reductive subalgebra $\fp$
and   a nilpotent ideal $\fn$;

  b) the action of $\fp$ on $\fn$ is
 completely reducible.}

 \smallskip

Obviously,  Lemma 1 implies the Ado theorem.
Indeed, $\fg$ admits a decomposition
$$\fg=\fp'\oplus (\fp''\ltimes\fn)$$
where $\fp'$, $\fp''$ are reductive
subalgebras
and the action of $\fp''$ on $\fn$
is faithfull. After this, it is sufficient
to apply the construction of p.1.



{\sc Remark.} The Ado theorem implies Lemma 1
modulo the Chevalley construction of algebraic
envelope of a Lie algebra.
But Lemma 1 itself can be easily
proved directly.

\smallskip

{\bf 3. Killing lemma.}
Let $\fg$ be a Lie algebra, let $d$ be its derivation.
For an eigenvalue $\lambda$, denote by
$\fg_\lambda$ its root subspace
$\fg_\lambda=\cup_k\ker(d-\lambda)^k$;
we have $\fg=\oplus \fg_\lambda$.
As it was observed by Killing,
$x\in \fg_\lambda$, $y\in \fg_\mu$ implies
$[x,y]\in \fg_{\lambda+\mu}$.

Thus the Lie algebra $\fg$ admits the gradation
by the eigenvalues of $d$.
Consider the gradation operator $d_s:\fg\to\fg$
defined by $d_s v=\lambda v$ if  $v\in\fg_\lambda$.
Obviously, $d_s$ is a derivation,
and $dd_s=d_sd$.
We also consider  the derivation $d_n:=d-d_s$,
this operator is nilpotent (the equality
$d=d_n+d_s$ is called
the Jordan--Chevalley decomposition).
 Clearly,
\begin{align}
\ker d_s\supset \ker d;\qquad\ker d_n\supset \ker d;\\
{\rm im}\,\, d_s\subset{\rm im}\,\, d_s;
\qquad {\rm im}\,\, d_n\subset{\rm im}\,\, d_s
\end{align}


{\bf 4. Elementary expansions.}
Let $\fq$ be a Lie algebra, let $I$ be an ideal
of codimension 1.
Let $x\notin I$.
Denote by $d$ the operator
$\Ad_x:I\to I$. Consider the corresponding
pair of derivations
$d_s$, $d_n$.
Consider the space
 $$\fq'=\C y + \C z + I$$
where $y$, $z$ are formal vectors.
We equip this space with a structure of
a Lie algebra
by the rule
$$
[y,z]=0, \qquad
[y,u]=d_su,\qquad [z,u]=d_nu,
\qquad \text{for all $u\in I$}
$$
and the commutator of $u,v\in I$
is the same as it was in $I$.

The subalgebra $\C(y+z)\oplus I\subset \fq'$ is
isomorphic $\fq$.
We say that $\fq'$ is an {\it elementary
expansion}
of $\fq$.

Obviously, $[\fq',\fq']=[\fq,\fq]$.









For a general Lie algebra,
 the required embedding to a semidirect
 product can be obtained by a sequence
of elementary expansions.


 \smallskip

 {\bf 5. Proof of Lemma 1.}
 Let $\fq$ be a Lie algebra. Let $\fh$ be its
 Levi part, and $\fr$ be the radical.
 Denote by $\fm$ the nilradical of
 $\fq$, i.e., $\fm=[\fq, \fr]$;
 recall that $\fm$ is a nilpotent ideal,
 and $[\fq,\fq]=\fh\ltimes\fm$
 (see \cite{Dix}, 1.4.9).

Consider a
  nilpotent ideal $\fn$ of $\fq$ containing
  the nilradical $\fm$.
 Consider a subalgebra $\fp\supset\fh$
 such that the adjoint
  action of $\fp$ on $\fq$ is 
  completely reducible
 and $\fp\cap \fn=0$;
 for instance, the can choice
 $\fn=\fm$, $\fp=\fh$.

Obviously,
 the $\fq$-module $\fq/(\fp\ltimes\fn)$
is trivial. Consider any
 subspace $I$ of codimension 1
  containing $\fp\ltimes\fn$,
  obviously $I$
is an ideal in $\fq$.
Since the action of $\fp$ on $\fq$
is completely reducible,
 there exists a $\fp$-invariant
complementary subspace
for $I$. Let $x$ be an element
of this subspace. Since the $\fp$-module
$\fq/I$ is trivial, $x$ commutes with
$\fp$.
We apply the elementary expansion to these
data.

We obtain the new algebra $\fq'=\C y + \C z+ I$
with the  nilpotent
ideal $\fn'=\C z + \fn$
and with the reductive subagebra $\fp'=\C y\oplus\fp$
(by (1), $y$ commutes with $\fp$).

It remains to notice that
$$
\dim \fq'-\dim \fp'-\dim \fn'=
\dim \fq-\dim \fp-\dim \fn-1
$$
and we can repeat the same construction.

\bigskip

{\sf Addres:

Math.Physics Group,
ITEP, B.Cheremushkinskaya, 25, Moscow, Russia

\&

Independent University of Moscow

\&

ESI, Vienna, Austria (January, 2002).}

\smallskip

neretin@main.mccme.rssi.ru

\end{document}